\documentclass[11pt, reqno]{amsart}
\usepackage{setspace}
\usepackage{textcomp}
\usepackage{amsmath,amssymb, mathrsfs}
\usepackage{bm}
\usepackage{fancyhdr}
\usepackage{lastpage}
\usepackage{color}
\usepackage{graphicx}
\usepackage{subfigure}
\usepackage{xcolor}
\usepackage[draft=false,colorlinks=true,CJKbookmarks=true,linkcolor=red,citecolor=blue,urlcolor=blue]{hyperref}
\usepackage{longtable}
\usepackage{amsthm}
\usepackage{lineno}
\usepackage{galois}

\usepackage{enumitem}
\usepackage{cite}

\newcommand{\N}{{\mathbb{N}}}
\newcommand{\Q}{{\mathbb{Q}}}
\newcommand{\R}{{\mathbb{R}}}

\renewcommand{\P}{{\mathbb{P}}}

\newcommand{\F}{{\mathcal{F}}}

\renewcommand{\H}{{\mathcal{H}}}
\newcommand{\G}{{\mathcal{G}}}
\newcommand{\C}{{\mathcal{C}}}

\DeclareMathAlphabet{\mathsfsl}{OT1}{cmss}{m}{sl}

\makeatletter
\newcommand{\neutralize}[1]{\expandafter\let\csname c@#1\endcsname\count@}
\makeatother

\newtheorem{thm}{Theorem}

\theoremstyle{plain}

\newtheorem{prop}[thm]{Proposition}
\newtheorem{cor}[thm]{Corollary}
\newtheorem{lem}[thm]{Lemma}
\newtheorem{pb}{Problem}
\newtheorem*{pb*}{Problem}

\theoremstyle{definition}
\newtheorem{defn}[thm]{Definition}

\theoremstyle{remark}
\newtheorem*{rmk}{Remark}

\begin{document}


\title[]{On Rainbow Thresholds}

\author{Jie Han }
\thanks{JH. School of Mathematics and Statistics and Center for Applied Mathematics, Beijing Institute of Technology, email: han.jie@bit.edu.cn}
\author{Xiaofan Yuan}
\thanks{Corresponding author: XY. School of Mathematical and Statistical Sciences, Arizona State University, email: xiaofan.yuan@asu.edu}
\thanks{JH is partially supported by the Natural Science Foundation of China (12371341).}

\begin{abstract}
Resolving a recent problem of Bell, Frieze, and Marbach, we establish both the threshold result of Frankston--Kahn--Narayanan--Park, and its strengthening by Spiro, in the rainbow setting.
This has applications to the thresholds for rainbow structures in random graphs where each edge is given a uniformly random color from a set of given colors.
\end{abstract}
\maketitle

\section{Introduction}

Threshold functions are central to the study of random discrete structures.
In the last few years, we have witnessed celebrated breakthroughs, e.g. by Frankston, Kahn, Narayanan, and Park~\cite{FKNP}, and Park and Pham~\cite{PP}.
The goal of this note is to extend this study to the rainbow setting, which has been the focus of a series of papers~\cite{BaFr, DEF, Bradshaw, FMc}, and formulated recently in~\cite{BFM}.

\subsection{Thresholds}
We start with some notation following~\cite{FKNP, T10}.
Given a finite set $X$, a family $\mathcal F\subseteq 2^X$ is called \textit{increasing} if $B \supseteq A \in \F \Rightarrow B \in \F$.
For a given $X$ and $p \in [0, 1]$, $\mu_p$ is the product measure on $2^X$ given by $\mu_p(S): = p^{|S|}(1-p)^{|X\setminus S|}$ for any $S\subseteq X$. 
For an increasing $\F$, the \textit{threshold} $p_c(\F)$ is the unique $p$ for which $\mu_p(\F) = \frac{1}{2}$. We say $\F$ is \textit{$p$-small} if there is a $\G \subseteq 2^X$ such that $\F \subseteq \langle \G \rangle := \{T : \exists S \in \G, S \subseteq T\}$ and
$\sum_{S \in \G} p^{|S|} \leq \frac{1}{2}$.
Then $q(\F) := \max\{p : \F$ is $p$-small$\}$, which is defined as the \textit{expectation-threshold of} $\F$.
Note that $q(\F)\le p_c(\mathcal F)$ by definition.
Let $\ell (\F)$ be the maximum size of minimal members of $\F$.

Talagrand~\cite{T10} introduced the following LP relaxation of ``$p$-smallness":
we say $\F$ is \textit{weakly $p$-small} if there is a function $g: 2^X \to \mathbb{R}^+$ such that
\[	
	\sum_{S\subseteq T} g(S) \geq 1  \,\, \forall T\in \F \quad \text{ and }\quad \sum_{S\subseteq X} g(S)p^{|S|} \leq \frac{1}{2}.\]
Then the \textit{fractional expectation-threshold} for $\F$, denoted by $q_f(\F)$, is defined as:
\[
q_f(\F) := \max\{p : \F \text{ is weakly } p\text{-small}\}.
\]


The following result was 
conjectured by Talagrand~\cite{T10}, and resolved by Frankston, Kahn, Narayanan, and Park~\cite{FKNP} recently, 
built on an ingenius approach on the Sunflower Lemma of Alweiss, Lovett, Wu, and Zhang \cite{ALWZ}.


\begin{thm}\cite{FKNP}\label{FKNPthm}
	There exists a constant $K$ such that for every finite $X$ and increasing $\mathcal F\subseteq 2^X$, 
 \[
p_c(\mathcal F) \le K q_f(\mathcal F) \log \ell(\mathcal F).
 \]
\end{thm}

This is improved by another recent result of Park and Pham~\cite{PP}, who showed a similar inequality resolving a conjecture of Kahn and Kalai~\cite{KK}, that 
\begin{equation}
 \label{eq:qpp}
     p_c(\mathcal F) \le K q(\mathcal F) \log \ell(\mathcal F).
 \end{equation}


\subsection{Rainbow Thresholds}
The goal of this note is to establish rainbow versions of these results, as recently asked by Bell, Frieze, and Marbach~\cite{BFM}, who also made a first attempt. 
A motivating question for thresholds in the rainbow setting could be the following (e.g. for random graphs).
\begin{pb}
    Determine the optimal $p=p(n)$ such that if the edges of $G(n,p)$ are randomly colored using $c\ge n$ colors, then the resulting graph contains a Hamilton cycle where each color is used at most once.
\end{pb}

In particular, Problem 1 has been resolved by Bal and Frieze~\cite{BaFr} while previous results have to use more colors (e.g., they proved it for $c\ge (1+o(1))n$ colors). 
See also~\cite{FeKr} for a more accurate bound on $p$.

To introduce this problem in the abstract setting we first define our (random) object and its threshold.
Given an increasing family $\mathcal F$, let $\H_\F$ be the family of minimal elements of $\mathcal F$.
We note that since $\F=\langle \H_\F\rangle$, $\F$ and $\H_\F$ are uniquely determined by each other.
Take an integer $k\ge \ell(\mathcal F)$ and let $[k]:=\{1,2,\dots, k\}$.
Our random object will be a randomly colored binomial random subset of $X$, that is, each $x\in X$ is given a uniformly random color in $[k]$ independently, and included in the random subset with probability $p$ independently of others (and independent of the color assignments). 
To describe the distribution of this process, we define a function $\mu_p^k: 2^{X'}\to [0,1]$, where $X':=X\times [k]$.
Given a set $S=\{(x_1, i_1), \dots, (x_t, i_t)\}\subseteq X'$, let
\[
\mu_p^k(S): = \begin{cases}
(p/k)^t (1-p)^{|X|-t}  & \text{ if all $x_1, \dots, x_t$ are distinct,} \\
0 & \text{ otherwise.}
\end{cases} 
\]
Note that $\mu_p^k(S)$ is equal to the probability that the outcome is exactly $S$,
if we color the elements of $X$ independently and uniformly at random with $[k]$ and take each $x\in X$ with probability $p$ independently.
Now if we extend this definition to subsets in $2^{X'}$ by $\mu_p^k(\F') := \sum_{S \in \F'} \mu_p^k(S) $ for $\F'\subseteq 2^{X'}$, then $\mu_p^k$ is a probability measure on $2^{X'}$ for $p\in (0,1)$.
{We also use $\mu'_p$ to denote this measure when not emphasizing the value of $k$.}

Given $\F\subseteq 2^X$, we win if this colored subset contains a \textit{rainbow edge} of $\H_\F$, that is, a subset with all vertices having distinct colors from $[k]$.
Let $\F^{\textrm{rb}}\subseteq 2^{X'}$ denote the increasing family generated by all the rainbow edges of $\H_\F$, and let the \textit{$k$-color rainbow threshold} $p_c^k(\mathcal F)$ for $\F$ be the value of $p$ such that $\mu'_p(\F^{\textrm{rb}}) = 1/2$ (where $k$ is the number of colors and $k\ge \ell(\F)$). 
We shall prove in Section 2 that for any non-trivial increasing family $\F'\subseteq 2^{X'}$, $\mu'_p(\F')$ is strictly increasing as $p$ increases in $(0,1)$, which guarantees that $p_c^k(\mathcal F)$ is well-defined.
Note that one may think of $p_c(\F)$ as $p_c^{\infty}(\F)$, namely, when $k$ tends to infinity, with probability tending to 1 all vertices of $X$ have distinct colors.

We indeed have 
\(
q(\F) \leq q_f(\F) \leq p_c(\F) \le   p_c^k(\F)\) for increasing $\F$ and $k\ge \ell (\F)$ (the second inequality was explained in~\cite{FKNP} and the third inequality will be shown in Section~\ref{sec:5}). 

Here is the main result of this note, which is a rainbow version of Theorem~\ref{FKNPthm}.

\begin{thm} \label{main}
There exists a constant $K'$ such that for any finite $X$, increasing family $\mathcal F\subseteq 2^X$, and integer $k\ge \ell(\mathcal F)$, we have
 \[
q_f(\F)\le p_c^k(\mathcal F) \le K' q_f(\mathcal F) \log \ell(\mathcal F).
 \]
\end{thm}

Note that by 
\eqref{eq:qpp},
we have $q_f(\F) \le p_c(\F)\le K q(\F)\log \ell(\F)$.
Therefore we get an immediate corollary 
\begin{equation}
\label{eq:log}
q(\F)\le p_c^k(\mathcal F) \le K'' q(\mathcal F) \log^2 \ell(\mathcal F)
\end{equation}
for some universal constant $K''$ and increasing $\mathcal F\subseteq 2^X$.

Thus, for increasing $\F$, its rainbow threshold and expectation-threshold (and thus also its threshold)
are up to a poly-logarithmic factor away from each other.
Following~\cite{FKNP}, we next state an equivalent version of Theorem~\ref{main} that is more handy for applications in the random graph setting.

\subsection{Spreadness}
Throughout this note, a \emph{hypergraph} on $X$ is a collection of subsets of $2^X$, with repeats allowed (unless stated otherwise).
A probability measure $\nu$ on $2^X$ is {\it $q$-spread} if 
\[	\nu (\langle S \rangle) \le q^{|S|}, \quad \forall S \subseteq X.\]
In particular, if $\H$ is a hypergraph on $X$, and the uniform measure $\nu$ on $\H$ is $\kappa^{-1}$-spread, then we say $\H$ is {\it $\kappa$-spread}.
We say $\H$ is {\it $r$-bounded} if $|e|\le r$ for all $e\in \H$ (i.e., $\ell (\langle \H \rangle) \le r$). 
We use $X_p$ to denote a subset of $X$ where each $x\in X$ is included independently in $X_p$ with probability $p$. 
Talagrand~\cite{T10} observed that for increasing $\F$, $q\ge q_f(\F)$ implies that there is a $(2q)$-spread measure on $\F$, and therefore, for applications, it usually \emph{suffices} to study/verify the spreadness condition.

Now we are ready to state an application-friendly version of our main result, phrased in the spreadness condition.
Note that its original (``uncolored'') version was proved in~\cite{FKNP}.

\begin{thm}\label{main_p}
There is an absolute constant $C>0$ such that the following holds.
Let $\H$ be an $r$-bounded $\kappa$-spread (multi-)hypergraph.
	Let $X=V(\H)$ and let the elements of $X$ be colored independently and uniformly at random from $[k]$, where $k\ge r$. 
	If 
	\[	p\ge \frac{C\log r}{\kappa},	\]
	then with probability $1-o_{r\to \infty}(1)$, $X_p$ contains a rainbow edge of $\H$. 
\end{thm}

Bell, Frieze, and Marbach~\cite{BFM} proved Theorem~\ref{main_p} with an additional assumption that $\kappa=\Omega(r)$ and conjectured that it could be removed.
Quick consequences of Theorem~\ref{main_p} include rainbow thresholds for Hamilton cycles in random (hyper)graphs and bounded degree spanning trees in random graphs. 
These results are already covered by the result of~\cite{BFM}, and the first one was established earlier by Bal and Frieze~\cite{BaFr} and Dudek, English and Frieze~\cite{DEF}.
We also provide new applications of Theorem~\ref{main_p} in Section~\ref{sec:appl}.

We also establish a similar result for the stronger spreadness, namely, the \textit{Spiro spread}~\cite{S23} (see also the \textit{superspread} by Espuny D\'iaz and Person~\cite{EP}).
This brings us immediate new applications, for example, the rainbow thresholds for $r$-th power of Hamilton cycles in random graphs, and the rainbow thresholds for Hamilton $\ell$-cycle in random hypergraphs for $\ell \ge 2$, again because the corresponding Spiro spread results in random hypergraphs have been established accordingly (see Section~\ref{sec:appl}).


\begin{defn}
	Let $0<q\le 1$ be a real number and $r_1>\cdots >r_\lambda$ be positive integers. We say that a hypergraph $\H$ on $V$ is \textit{$(q;r_1,\ldots,r_\lambda)$-spread} if $\H$ is non-empty, $r_1$-bounded, and if for all $A\subseteq V$ with $d(A)>0$ and $r_{i}\ge|A|\ge r_{i+1}$ for some $1\le i< \lambda$, we have for all $j\ge r_{i+1}$ that
	\[M_j(A):=|\{S\in \H:|A\cap S|= j\}|\le q^j |\H|.\]
\end{defn}

We prove the Spiro-spreadness version of Theorem~\ref{main_p} as follows, and if we ignore the colors then it is exactly the result obtained by Spiro~\cite{S23} (see Theorem~\ref{Spiromain}).

\begin{thm}
\label{main_p_spiro}
There is an absolute constant $K>0$ such that the following holds.
Let $\H$ be an $r_1$-uniform $(1/\kappa;r_1,\ldots,r_\lambda)$-spread (multi-)hypergraph with $r_\lambda=1$.
	Let $X=V(\H)$ and let the elements of $X$ be colored independently and uniformly at random from $[k]$, where $k\ge r_1$. 
	If 
	\[	p\ge \frac{C\lambda}{\kappa}	\]
        with $C\ge K$,
	then with probability $1-K/(C\lambda)$, $X_p$ contains a rainbow edge of $\H$. 
\end{thm}

\subsection{Applications}
Now we give two applications of Theorem~\ref{main_p_spiro}.

The first one is on Hamilton $\ell$-cycles in random hypergraphs.
A \textit{$k$-uniform $\ell$-cycle of length $t$} is a $k$-uniform hypergraph with $t(k-\ell)$ vertices and $t$ edges, where the vertex admits a cyclic order $v_1,\dots,v_{t(k-\ell)}$ and the edges are $\{v_{i(k-\ell)+1},\dots,v_{i(k-\ell)+k}\}$ for $i\in [t]$ (the addition on the indices are modulo $t(k-\ell)$).
When $\ell \ge 2$, it is known that the threshold for Hamilton $\ell$-cycles in binomial random $k$-uniform hypergraphs is $n^{-(k-\ell)}$.
With a Spiro spread condition proved by Kelly, M\"uyesser, and Pokrovskiy \cite{KMP}, Theorem~\ref{main_p_spiro} implies that the rainbow threshold for Hamilton $\ell$-cycles in random $k$-uniform hypergraphs is also $n^{-(k-\ell)}$.

\begin{thm}\label{ellcycle}
For all $k,\ell \in \mathbb N$ with $2\le \ell <k$ and $\varepsilon >0$, there exists $K'=K'(k, \varepsilon)$ such that for all sufficiently large $n$ with $(k-\ell)|n$, 
if the edges of $G^{(k)}(n,p)$ are colored independently and uniformly at random with $q\ge n/(k-\ell)$ colors, where $p=K'n^{\ell-k}$, then
the resulting hypergraph
contains a rainbow Hamilton $\ell$-cycle.
\end{thm}

The second application is on the $r$-th power of Hamilton cycles.
Given a graph $G$, the $r$-th power of $G$ is obtained from $G$ by adding edges to all pairs of vertices with distance at most $r$.
Very recently Bell and Frieze~\cite{BF} studied the rainbow threshold for $r$-th power of Hamilton cycles in random graphs, for $r\ge 2$.
They showed that $n^{-1/r}$ is the threshold if the number of colors is at least $(1+\varepsilon)rn$.
We use Theorem~\ref{main_p_spiro} to reduce the number of colors to $rn$ (best possible).

\begin{thm}\label{powerofham}
For all $\varepsilon >0$ and $r\ge 2$, there exists a constant $K'=K'(r, \varepsilon)$ such that 
if the edges of $G(n,p)$ are colored independently and uniformly at random with $q\ge rn$ colors, where $p=K'n^{-1/r}$, then
with probability at least $1-\varepsilon$, 
the resulting graph contains a rainbow $r$-th power of Hamilton cycle.
\end{thm}

In general, Theorem~\ref{main_p_spiro} can determine the rainbow threshold for (spanning) subgraph $F$ if one can establish the corresponding Spiro-spreadness for the hypergraph formed by all copies of $F$ in the complete (hyper)graph.
Theorems~\ref{ellcycle} and~\ref{powerofham} are proved in Section~\ref{sec:appl}, where we also present two applications of Theorem 3 on randomly colored random sparsifications of dense (hyper)graphs.


\subsection{Transversal versions}
By the above observation of Talagrand and routine computations to push the error probabilities (see~\cite{FKNP}), Theorem~\ref{main_p} implies the upper bound in Theorem~\ref{main}.
To prove Theorem~\ref{main_p}, we first prove the following transversal version and then couple these two models following 
an ingenious coupling idea of McDiarmid \cite{MD80}, see other applications by Ferber \cite{A18} and Ferber--Krivelevich~\cite{FeKr}.
The transversal variant was already mentioned in~\cite{BFM} as a possible route (see Remark 2 therein).

\begin{thm} \label{trans2}
There is an absolute constant $C>0$ such that the following holds.
	Let $\H$ be an $r$-bounded, $\kappa$-spread (multi-)hypergraph. Let $X=V(\H)$ and $k\ge r$, and let $X' = X\times [k]$.
 If
	\[	p\ge \frac{C\log r}{ \kappa},	\]
	then with probability $1-o_{r\to \infty}(1)$, $X'_{p/k}$ contains a set $S'
	= \{(x_1, i_1), (x_2,i_2), \dots, (x_t, i_t)\}$, where 
	$\{x_1, x_2, \dots, x_t\}$ is an edge of $\H$ and $i_1, i_2,\dots, i_t$ are distinct numbers in $[k]$.
\end{thm}

Theorem~\ref{trans2} is called a \emph{transversal} version because we can naturally consider a family  of hypergraphs $\mathscr F = \{\F_1, \F_2, \dots, \F_k\}$, where $V(\F_i) = X\times \{i\}$, and $\{(x_1, i), (x_2,i), \dots, (x_t, i)\}\in \F_i$ if and only if $\{x_1, x_2, \dots, x_t\}\in \F$. 
Then our target objects are the sets $\{(x_1, i_1),\dots, (x_t, i_t)\}$ such that $\{x_1,\dots,x_t\}\in \mathcal F$ and $i_1,\dots, i_t$ are all distinct. 
 In particular, each of them contains at most one vertex from $V(\F_i)$ for each $i$.
In the special case $k=r$ and $\H$ is $r$-uniform, the target objects are the transversals in $\mathscr F$, that is, each target object contains exactly one vertex from $V(\F_i)$ for each $i$.

\medskip
The rest of this paper is organized as follows.
We prove Theorem~\ref{trans2} in Section \ref{sec:3}, and Theorem~\ref{main_p} and Theorem~\ref{main_p_spiro} in Section \ref{sec:4}. The proof of Theorem~\ref{main} is given in Section \ref{sec:5}.
We give some applications of our results in Section \ref{sec:appl} and conclude with some further directions in Section \ref{sec:rem}.
In the appendix we show that $\mu'_p$ is strictly increasing for increasing families when $p\in (0,1)$, justifying our definition of the rainbow threshold.

\section{Proof of Theorem~\ref{trans2}}
\label{sec:3}

Write $(k)_r:=k (k-1) \cdots (k-r+1)$ for $k, r\in \mathbb N$. 
The proof of Theorem~\ref{trans2} uses the fact that the rainbow colorings of a hyperedge $E$ are evenly distributed, that is, there are exactly $(k)_{|E|}$ such colorings given $k$ colors to use.
We use this to reduce the problem to its uncolored version.
Note that a similar idea was observed (and used) in~\cite{BFM}.

We need the following result in \cite{FKNP}.

\begin{thm}\cite[Theorem 1.6]{FKNP}\label{FKNP1.6}
There is an absolute constant $K>0$ such that 
for any $r$-bounded, $\kappa$-spread multi-hypergraph $\mathcal{J}$ on $X$, 
a uniformly random $((K \kappa^{-1}\log r)|X|)$-element subset of $X$ belongs to $\langle \mathcal{J} \rangle$ with probability $1-o_{r \to \infty}(1)$.
\end{thm}

For an edge $E$ in a multi-hypergraph $\H$, let $m_\H(E)$ be its multiplicity, namely, the number of times it appears in $\H$.

 \begin{proof}[Proof of Theorem~\ref{trans2}]
	Let $\H'$ be the hypergraph on $X'$ where 
	$S' = \{(x_1, i_1), (x_2,i_2), \dots, (x_t, i_t)\}$ is an edge in $\H'$ if and only if 
	$\{x_1, x_2, \dots, x_t\}$ is an edge in $\H$ and $i_1, i_2,\dots, i_t$ are distinct numbers in $[k]$, and the multiplicity of $S'$ in $\H'$ is equal to the multiplicity of $\{x_1, x_2, \dots, x_t\}$ in $\H$, that is, $m_{\H'}(S')=m_{\H}(\{x_1, x_2, \dots, x_t\})$.
	We define an auxiliary multi-hypergraph $\H''$ by the following: for each $S'$ in $\H'$, we include $(k-|S'|)_{r-t}$ copies of $S'$ in $\H''$. That is, $\H''$ and $\H'$ have the same set of members (edges), and for every $S'\in \H''$, $m_{\H''}(S')=(k-|S'|)_{r-t}m_{\H'}(S')$.
	
	Note that for each edge $\{x_1, x_2, \dots, x_t\}$ in $\H$, there are $(k)_{t}$ choices for $i_1, i_2,\dots, i_t$ such that $S'= \{(x_1, i_1), (x_2,i_2), \dots, (x_t, i_t)\}$ is an edge in $\H'$.
 Thus each edge $\{x_1, x_2, \dots, x_t\}$ in $\H$ corresponds to $(k)_t\cdot (k-t)_{r-t}= (k)_r$
	edges in $\H''$. Hence we have $|\H''| = (k)_r |\H|$.
	
Now we show that $\mathcal H''$ is $(k\kappa/e)$-spread. Indeed, let $S=\{(x_1, i_1),\dots, (x_s, i_s)\}$ be a subset of $X'$. 
We may assume that $x_1, x_2, \dots, x_s$ are distinct in $X$ and $i_1, i_2,\dots, i_s$ are distinct in $[k]$, as otherwise $|\H'' \cap \langle S \rangle|=0$. 	
	Let 
	\[	S\subseteq T = \{(x_1, i_1), (x_2, i_2), \dots, (x_s, i_s), (x_{s+1}, i_{s+1}), \dots, (x_{t}, i_t)\}\subseteq X'	\]
	and note that $T\in \H''$ if and only if $T^* = 
	\{x_1, x_2, \dots, x_t\}$ is an edge of $\H$ and $ \{i_1, i_2,\dots, i_t\}$ consists of $t$ distinct numbers in $[k]$. 
So, for each edge $T^*$ in $\H$, there are $(k-s)_{t-s}$ choices for $\{i_{s+1}, \dots, i_t\}$, each of which corresponds to $(k-t)_{r-t}$ (same) edges.
Then each $T^* \in \H\cap \langle \{x_1, x_2, \dots, x_s\} \rangle$ corresponds to $(k-s)_{t-s}\cdot (k-t)_{r-t} = (k-s)_{r-s}$
	edges in $\H''\cap \langle S \rangle$.
Thus we have
    \begin{align}\label{extspread}
		|\H''\cap \langle S \rangle|  = \sum_{\{x_1, x_2, \dots, x_s\} \subseteq T^*\in \H} (k-s)_{r-s} \le \frac{|\H|}{\kappa^{s}} (k-s)_{r-s} = \frac{|\H''|}{\kappa^{s}(k)_s} \le \frac{e^s|\H''|}{(k \kappa)^{s}},
    \end{align}
  where we used that $(k)_s \ge (k/e)^s$.
	Hence, $\H''$ is $(k\kappa/e)$-spread.
	
	Let $K$ be the constant from Theorem~\ref{FKNP1.6}.
	Let $C = 2eK$
	and let 
	\[ p\ge \frac{C\log r}{\kappa} = \frac{(C/e)\log r}{\kappa/e} = \frac{2K\log r}{\kappa/e}. \]
Applying Theorem~\ref{FKNP1.6} with $\H''$ in the place of $\mathcal{J}$ implies that
a uniformly random $(p|X'|/2k)$-element subset of $X'$ contains an edge in $\H''$ with probability $1-o_{r \to \infty}(1)$.
Standard concentration arguments give that $|X_p'|\ge p|X'|/2k$ holds with probability $1-o_{|X|\to \infty}(1)$.
Conditioning on this, we can take a random subset of size exactly $p|X'|/2k$.
Since when conditioning on the size of the outcome, the binomial distribution reduces to the hypergeometric distribution, 
we obtain that  
with probability $1-o_{r \to \infty}(1)$, $X'_{p}$ contains a desired edge in $\H''$ (and hence in $\H'$ as well).
 \end{proof}

\begin{rmk}\label{spiro}
The same approach applies to Spiro spread properties. We prove the following as a Spiro spread analoge of Theorem~\ref{trans2}.
\end{rmk}

\begin{thm} \label{trans2_spiro}
There is an absolute constant $K>0$ such that the following holds.
	Let $\kappa \ge 6$ and $\H$ be an $r_1$-uniform $(1/\kappa;r_1,\ldots,r_\lambda)$-spread (multi-)hypergraph with $r_\lambda=1$. Let $X=V(\H)$ and $k\ge r$, and let $X' = X\times [k]$.
 If
	\(	p\ge C\lambda/ {\kappa}	\)
	with $C\ge K$, then with probability $1-K/(C\lambda)$, $X'_{p/k}$ contains a set $S'
	= \{(x_1, i_1), (x_2,i_2), \dots, (x_t, i_t)\}$, where 
	$\{x_1, x_2, \dots, x_t\}$ is an edge of $\H$ and $i_1, i_2,\dots, i_t$ are distinct numbers in $[k]$.
\end{thm}

\begin{proof}[Sketch of proof of Theorem~\ref{trans2_spiro}]
The proof follows that of Theorem~\ref{trans2}.
We define $\H', \H''$ in the same way.
(As $\H$ is uniform, we indeed have $\H''=\H'$.)
Then we have $|\H''|=(k)_r|\H|$.

To show that $\H''$ is $({k\kappa}/({3e});r_1,\dots,r_\lambda)$-spread, we follow the same proof but replace the inequality \eqref{extspread} by the following.
For $S\subseteq \mathcal H'$, let us write $S_X$ for the projection of $S$ onto $X$.
Now for all $A\subseteq X'$ with $r_{i}\ge|A|\ge r_{i+1}$ for some $1\le i< \lambda$, and for all $s\ge r_{i+1}$, we have
\begin{align*}
    |\{T\in \H'': |A\cap T|= s\}| 
    & = 
    \sum_{t=s}^{|A|}	
    |\{T\in \H'': |A\cap T| = s,~|T_X\cap A_X| = t\}| 
    , 
    \\
    & \le 
    \sum_{t=s}^{|A|} 
    \left(\sum_{T_X\in \H: ~|T_X\cap A_X| = t,~|A\cap T|=s} \binom{t}{s} (k-s)_{r-s} \right)
    \\
    & \le (k-s)_{r-s}\left(\sum_{t=s}^{|A|} |\{T_X\in \H: |A_X\cap T_X|=t\}|  \cdot \binom{t}{s} \right)\\
    & \le (k-s)_{r-s}\left(\sum_{t=s}^{|A|} \frac{|\H|}{\kappa^{t}} \cdot 2^t \right)  \\
    & \le \frac{|\H''|}{(k)_s}\cdot \left(\frac{2}{\kappa}\right)^s\frac{1}{1-2/\kappa} \le \frac{|\H''|}{(k)_s}\left(\frac{3}{\kappa}\right)^s \le \frac{|\H''|}{(k \kappa/3e)^{s}}
    \end{align*}
where we used $\kappa\ge 6$ and $(k)_s \ge (k/e)^s$, and in the second line, we used that given $T_X\in \H$, the number of choices for different $T$ corresponding to $T_X$ is at most $\binom{|A_X\cap T_X|}{s} (k-s)_{r-s} = \binom ts (k-s)_{r-s}$, that is, we first choose $A\cap T$ (by choosing $(A\cap T)_X\subseteq A_X\cap T_X$, a set of $s$ vertices) and then choose $T$ (color the other $r-s$ vertices). 

We then continue the rest of the proof, except that apply the following result of Spiro~\cite{S23} (Theorem~\ref{Spiromain}) instead of Theorem~\ref{FKNP1.6}.
\end{proof}

\begin{thm}[\cite{S23}, Theorem~1.3]
\label{Spiromain}
There exists an absolute constant $K_0$ such that the following holds. Let $\H$ be an $r_1$-uniform $(q;r_1,\ldots,r_\lambda)$-spread hypergraph on $V$ with $r_\lambda=1$.  If $W$ is a set of size $C\lambda q|V|$ chosen uniformly at random from $V$ with $C\ge K_0$, then
	\[\Pr[W\text{ contains an edge of }\H]\ge 1-\frac{K_0}{C\lambda}.\]

\end{thm}

\section{Proof of Theorems~\ref{main_p} and~\ref{main_p_spiro}}
\label{sec:4}

We prove the following result that couples our two models, the randomly color model and the transversal model, following an ingenious idea of McDiarmid~\cite{MD80}.

\begin{lem} \label{couple}
	Let $X$ be a finite set and let $\H$ be an $r$-bounded hypergraph on $X$. For $k\ge r$, let $X' = X\times [k]$. 
	Let $\H'$ denote the family consisting of all possible rainbow copies of $\H$ using the color set $[k]$.
	Let $p\in (0,1)$ and $p'=p/k$.
	Consider the following two events:
	\begin{itemize}
		\item Event T: $X'_{p'}$ contains an edge of $\H'$.
		\item Event C: Given a random coloring on $X$ where the elements of $X$ are independently and uniformly colored from a set $\C = [k]$. Then $X_p$ contains a rainbow edge of $\H$.
	\end{itemize}
	Then $\P[T] \le \P[C]$.
\end{lem}

\begin{proof}
We define intermediate random samplings ``between'' those two events. 
Let $n = |X|$, and we may assume $X=[n]$. 
We define the random sets $X^0, X^1, \dots, X^n$ by the following: \\
In $X^i$, for each $j\le i$, we include $(j, s)$ with probability $p$, where $s$ is a uniformly chosen color in $[k]$
for $j$ in Event C; 
for each $j>i$, we include $(j,s)$ in $X^i$ with probability $p'$ independent for all colors $s$ in $[k]$. 
Then, in particular, $X^0 = X'_{p'}$ as in Event T, and $X^n$ contains an edge of $\H'$ if and only if Event C occurs.

Now it suffices to show
\[
	\P[X^i \text{ contains an edge of } \H'] \ge \P[X^{i-1} \text{ contains an edge of } \H']
	\]
for all $i\in [n]$.

Note that $X^i$ and $X^{i-1}$ only differ at the pairs with $i$ as the first coordinate. 
We can divide into the following three cases:
\begin{itemize}
		\item [a.] $X^{i-1}\setminus \{(i,s): s\in [k]\}$ contains an edge of $\H'$ (i.e., there exists such an edge not using vertex $i$)
		\item [b.] $X^{i-1}\cup \{(i,s): s\in [k]\}$ does not contain an edge of $\H'$ (i.e., there does not exist such an edge of $\H'$ even with all colors for vertex $i$ available.)
		\item [c.] Not in the case of (a) or (b). That is, $X^{i-1}$ contains an edge of $\H'$ or not depending on the occurrence of the second coordinate (the color) associated with vertex $i$.
	\end{itemize}
We consider conditional probabilities. Note that in case (a), $X^{i-1}$ contains an edge of $\H'$, and that edge also shows up in $X^{i}$. Thus 
\[
	\P[X^i \text{ contains an edge of } \H'\ | \ case (a)] = \P[X^{i-1} \text{ contains an edge of } \H'\ | \ case (a)] = 1.
	\]
Similarly, 
\[
	\P[X^i \text{ contains an edge of } \H'\ | \ case (b)] = \P[X^{i-1} \text{ contains an edge of } \H'\ | \ case (b)] = 0.
	\]
In case (c), for an arbitrary instance of $X^{i-1}\setminus \{(i,s): s\in [k]\}$, say $\tilde{X}$, we define 
\[
	D= D_{\tilde{X}}: = \left\{c\in [k]: \tilde{X}\cup \{(i, c)\} \text{ contains an edge of } \H'	\right\},
	\]
	that is, the set of colors for $i$ that completes an edge of $\H'$. 
 Note that $D\neq \emptyset$ as $k\ge r$ and we are not in case (b).
Then 
\[
	\P\left[X^{i-1} \text{ contains an edge of } \H'\ \big| 
	\ \tilde{X}\right] = 1- (1-p')^{|D|}, 
	\]
	that is, at least one color in $D$ occurs for $i$. 
On the other hand, 
\[
	\P\left[X^{i} \text{ contains an edge of } \H'\ \big| 
	\ \tilde{X}\right] = p \cdot \frac{|D|}{k},
	\]
	that is, the assigned color for $i$ is in $D$.
Note that $p'=p/k$, $|D|\le k$, and $p< 1$ imply
\begin{align*}
	p \cdot \frac{|D|}{k} - \left( 1- (1-p')^{|D|} \right) & = p \cdot \frac{|D|}{k} - 1+ (1-p/k)^{|D|}  \\
	& = \sum_{j=1}^{\infty} \left((p/k)^{2j}\binom{|D|}{2j} - (p/k)^{2j+1}\binom{|D|}{2j+1} \right).
	\end{align*}
Since  $p|D|/k< 1$, we have
	\[
		(p/k)^{2j}\binom{|D|}{2j} - (p/k)^{2j+1}\binom{|D|}{2j+1}\ge (p/k)^{2j}\binom{|D|}{2j} \left(1-\frac{p|D|}{k(2j+1)}\right) > 0.	\]
Thus we have for any instance $\tilde{X}$ in case (c),
\[
	 \P\left[X^{i} \text{ contains an edge of } \H'\ \big| 
	\ \tilde{X}\right] \ge \  \P\left[X^{i-1} \text{ contains an edge of } \H'\ \big| 
	\ \tilde{X}\right],
\]
giving that 
\[
	\P[X^i \text{ contains an edge of } \H'\ | \ case (c)] \ge \P[X^{i-1} \text{ contains an edge of } \H'\ | \ case (c)] .
	\]
Therefore, we have
\[
	\P[X^i \text{ contains an edge of } \H'] \ge \P[X^{i-1} \text{ contains an edge of } \H'],
	\]
and inductively on $i$, we conclude that $\P[T] \le \P[C]$.
\end{proof}

Now we are ready to prove Theorems~\ref{main_p} and~\ref{main_p_spiro}.

\begin{proof}
[Proof of Theorems~\ref{main_p} and~\ref{main_p_spiro}]
For Theorem~\ref{main_p}, take $C$ from Theorem~\ref{trans2}.
Let $p\ge C\log r / \kappa$ and $p'=p/k$.
Theorem~\ref{trans2} 
says that $X_{p'}'$ contains an edge of $\H$ with probability $1-o_{r\to \infty}(1)$.
Hence by Lemma~\ref{couple}, if we color $X$ uniformly at random with $k$ colors and take each $x\in X$ with probability $p$ independently, then with probability $1-o_{r\to \infty}(1)$ the outcome contains a rainbow edge of $\H$.

The same argument with Theorem~\ref{trans2_spiro} in place of Theorem~\ref{trans2} proves Theorem~\ref{main_p_spiro}.
\end{proof}

\section{Proof of Theorem~\ref{main}}
\label{sec:5}

Recall that the second inequality is given by Theorem~\ref{main_p} combined with the observation of Talagrand.
So it remains to prove the first inequality.

Note that it suffices to show $q_f(\F)\le p_c(\F) \le p_c^k (\F)$ where the first inequality is known. 
Let $\F^{all}\subseteq 2^{X'}$ denote the collection of all items in $\F$ with every possible coloring from $[k]$, that is, 
\[
\F^{all} := \left\{ \{(x_1, i_1), \dots (x_s, i_s) \}: \{x_1, \dots, x_s\} \in \F \text{ and } i_j\in [k] ~\forall j\in [s] \right\}.
\] Then we have $\mu_p(\F) = \mu_p^k(\F^{all})$ by definition, and $\mu_p^k(\F^{all}) \ge \mu_p^k(\F^{rb})$ as $\F^{all}\supseteq \F^{rb}$, for all $p\in (0,1)$ and positive integer $k$. 
Recall that for any non-trivial increasing family $\F\subseteq 2^{X}$, when $p$ increases in $(0,1)$, both $\mu_p(\F)$ and  $\mu_p^k(\F^{rb})$ are strictly increasing (the former one is observed in~\cite{FKNP} and the latter one is Lemma~\ref{mono}).
Thus, $p_c(\F) \le p_c^k (\F)$ by definition.
\qed


\section{Applications of Theorem~\ref{main_p} and Theorem~\ref{main_p_spiro}}
\label{sec:appl}

We first prove Theorems~\ref{ellcycle} and~\ref{powerofham}.
\subsection{Spiro-spread -- proofs of Theorems~\ref{ellcycle} and~\ref{powerofham}}


Recently, Kelly, M\"uyesser, and Pokrovskiy \cite{KMP} proved that for $\ell\ge 2$ the uniform measure on the family of Hamilton $\ell$-cycles in $K_n^{(k)}$ satisfies Spiro spread condition.
The following is a consequence of \cite[Proposition 6.4]{KMP}.

\begin{prop}\cite{KMP}
    \label{prop:ell-cycle-spreadness}
    For every 
    $k \in \mathbb N$, there exists $C=C(k)>0$ such that 
    for all $\ell \in \{2, \dots, k-1\}$ and all sufficiently large $n$ with $(k-\ell)|n$,
    the collection of all Hamilton $\ell$-cycles in $K_n^{(k)}$ forms an $(n/(k-\ell))$-uniform 
    $(C/n^{k-\ell}; n/(k-\ell), 1)$-spread hypergraph, 
    where $K_n^{(k)}$ denotes the complete $k$-uniform hypergraph on $n$ vertices. 
\end{prop}
With Theorem~\ref{main_p_spiro}, it implies the rainbow threshold for Hamilton $\ell$-cycles in random $k$-uniform hypergraphs is $n^{-(k-\ell)}$ as in Theorem~\ref{ellcycle}.



Secondly, for square of Hamilton cycles (the case of $r=2$), as pointed out in \cite{S23}, Kahn, Narayanan, and Park \cite{Hamiltonian} implicitly proved that the hypergraph ${\H}$ encoding squares of Hamiltonian cycles is a $(2n)$-uniform $(Cn^{-1/2};2n,C_0n^{1/2},1)$-spread hypergraph for some appropriate constants $C,C_0$. Using this and Theorem~\ref{main_p_spiro}, we obtain that $n^{-1/2}$ is the threshold for $2n$ colors.

For $r$-th power of Hamilton cycles, Chen, Li, Zhan, and the first author~\cite{CHLZ-threshold} recently proved the following result.
We note that it is also implicitly proved in~\cite{Joos-Lang-Sanhueza}.

\begin{thm}\cite{CHLZ-threshold}
Given $r\ge 2$, there exist constants $C',C_0'$ such that the collection of all $r$-th power of Hamilton cycles in $K_n$ forms an $(rn)$-uniform $(C'n^{-1/r};rn, n/2, C_0'n^{1-1/r},1)$-spread hypergraph.
\end{thm}

This together with Theorem~\ref{main_p_spiro} implies Theorem~\ref{powerofham}.

\subsection{Spread -- applications of Theorem~\ref{main_p}}
Below we collect some quick implications by combining Theorem~\ref{main_p} and some known spreadness computations.
We note that for these applications the result of Bell--Frieze--Marbach~\cite{BFM} (a weak version of Theorem~\ref{main_p}) is enough.

The following remark is useful in our applications.
 
 \begin{rmk}
Suppose $\nu$ is a $\kappa^{-1}$-spread measure on $2^X$ supported on $\H$.
Since $\Q$ is dense in $\R$, we may assume that $\nu$ takes values in $\Q$ (up to an arbitrarily small error term in the spread value, which could be taken care of by adjusting the constants). Now we remark that the assumption $\H$ is $\kappa$-spread can be replaced by that there exists a $\kappa^{-1}$-spread probability measure (distribution) on $\H$, by taking multiple edges in $\H$ where $\nu$ corresponds to the uniform distribution on the new edge set. 
\end{rmk}

\subsubsection{Spanning trees}

The authors of \cite{PSSS} proved 
the following result for bounded degree spanning trees.
Here $\delta$ denotes the minimum degree and $\Delta$ denotes the maximum degree of a graph.
\begin{lem}[Lemmas 7.2 \& 7.3 in \cite{PSSS}] \label{tree}
	For every $\Delta \in \N$ and $\delta > 0$ there exists $C_0=C_0(\Delta, \delta)>0$ such that the following holds for sufficiently large integer $n$.
 For every $n$-vertex graph $G$ with $\delta (G)\ge (1/2+\delta)n$ and every tree $T$ on $n$ vertices with $\Delta(T)\le \Delta$ there exists a $(C_0/n)$-spread distribution on graph embeddings of $T$ into $G$.
\end{lem}

Combining the above result with Theorem~\ref{main_p}, one can derive the following result on the threshold of a randomly colored binomial random subgraph of a given graph with a large minimum degree to contain a rainbow copy of a given bounded-degree spanning tree.
\begin{thm}\label{treethm}
	For every $\Delta \in \N$ and $\delta > 0$ there exists $C=C(\Delta, \delta)>0$ such that the following holds for sufficiently large integer $n$. Suppose that $G$ is an $n$-vertex graph satisfying $\delta(G)\ge (1/2+\delta)n$ and $T$ is an $n$-vertex tree with $\Delta(T)\le \Delta$. Suppose the edges of $G$ are colored independently and uniformly at random with $q\ge n-1$ colors. 
Then with probability $1-o_{n\to \infty}(1)$, $G_{C \log n/n}$ contains a rainbow copy of $T$, where $G_p$ denotes the spanning subgraph of $G$ where each edge of $G$ is retained with probability $p$. 
\end{thm}

\begin{proof}
	We may assume $n$ is sufficiently large. Then by Lemma~\ref{tree}, there exists a $(C_0/n)$-spread distribution on graph embeddings of $T$ into $G$.
Let $\H$ be a hypergraph whose vertices are the edges of $G$ and edges are the graph embeddings of $T$ into $G$. Then $\H$ is $(n-1)$-uniform and has a $(C_0/n)$-spread distribution. By Theorem~\ref{main_p}, there exists a constant $C_1$ such that if 
\[
p\ge  {C_1\log (n-1)} \frac{C_0}n ,
\]
then with probability $1-o_{n\to \infty}(1)$, $X_p$ contains a rainbow edge of $\H$.

Let $C= C_1C_0$, and then $C \log n/n > {C_1\log (n-1)}\frac{C_0}n $. Therefore, with probability $1-o_{n\to \infty}(1)$, $G_{C \log n/n}$ contains a rainbow copy of $T$. 
\end{proof}

\subsubsection{Matching in hypergraphs}
All hypergraphs considered in this subsection are simple, i.e., no repeats allowed.
Let $\H$ be a $k$-uniform hypergraph with vertex set $V$.  For any $T\subseteq V$, we use $d_{\H}(T)$ to denote the {\it degree}
of $T$ in $\H$, i.e., the number of edges of $\H$
containing $T$. For a positive integer $\ell$ such that $1\le \ell < k$, define 
$\delta_{\ell}(\H):= \min\left\{d_{\H}(T): T\in \binom{V}{\ell}\right\}$
to be the minimum {\it $\ell$-degree} of $\H$. 

For integers $\ell,k,n$ satisfying $1\le \ell< k$ and $n\in k\N$, let $t(n, k, \ell)$ be the smallest $d$ such that every $n$-vertex $k$-uniform hypergraph with $\delta_{\ell}(\H)\ge d$ contains a perfect matching. Define the {\it $\ell$-degree (Dirac) threshold} for perfect matchings in $k$-uniform hypergraphs to be
\[	\delta_{\ell,k}^+ := \lim_{k|n, n\to \infty} \frac{t(n,k,\ell)}{\binom{n}{k-\ell}}. \]

In the independent works of \cite[Theorem 1.5]{KKOP} and~\cite{PSSS}, the authors showed that there is a spread distribution on the perfect matchings of the $k$-uniform hypergraphs satisfying the Dirac condition. 
\begin{lem}\cite{KKOP, PSSS} \label{PMspread}
	Let $k,\ell$ be integers such that $1\le \ell< k$ and let $\varepsilon >0$. There exists $C_0=C_0(\ell,k,\varepsilon)$ such that the following holds for sufficiently large integer $n$. Let $\H$ be an $n$-vertex $k$-uniform hypergraph such that $k|n$ and $\delta_{\ell}(\H)\ge (\delta_{\ell,k}^+ +\varepsilon)\binom{n}{k-\ell}$.
	Then there exists a probability measure on the set of perfect matchings in $\H$ which is $(C_0/n^{k-1})$-spread.
\end{lem}

With Theorem~\ref{main_p}, we have the following corollary. The proof follows the lines in the above subsection and is omitted.

\begin{cor}
\label{coro}
	Let $k,\ell$ be integers such that $1\le \ell< k$ and let $\varepsilon >0$. There exists $C=C(\ell,k,\varepsilon)$ such that the following holds for sufficiently large integer $n$. Let $\H$ be an $n$-vertex $k$-uniform hypergraph such that $n\in k\N$ and $\delta_{\ell}(\H)\ge (\delta_{\ell,k}^+ +\varepsilon)\binom{n}{k-\ell}$. Suppose the edges of $\H$ are colored independently and uniformly at random with $q$ colors, where $q\ge n/k$. Then with probability $1-o_{n\to \infty}(1)$, $\H_{C\log n / n^{k-1}}$ contains a rainbow perfect matching.
\end{cor}

The sharp minimum $(k-1)$-degree condition for the above results is proven in \cite[Theorem 1.6]{KKOP}. Together with Theorem~\ref{main_p} this also gives the sharp minimum $(k-1)$-degree condition for the robustness of rainbow perfect matchings. Here we omit the very similar statements.





\section{Concluding Remarks}
\label{sec:rem}

There are at least three possible further problems left after this paper.
First, it is natural to try to prove a rainbow version of~\eqref{eq:qpp},
that is, to remove the extra logarithmic factor in~\eqref{eq:log}.
Following the proof idea of this note, it suffices to prove a ``transversal version", which we did not manage to do.
However, it would also follow from a conjecture of Talagrand~\cite[Problem 6.3]{T10}, who suggested that $q(\F)\ge q_f(\F)/K$ for some absolute constant $K$ and every increasing $\F$.

Second, it is not clear to us how to obtain tight(er) results on rainbow structures when the term $\log n$ can be reduced (but not \textit{dropped}) in the uncolored version.
For example, for the $K_r$-factor case, the uncolored version can be resolved (see~\cite[Section 7]{FKNP}) by using a nice coupling result of Riordan~\cite{Riordan} and converting the problem to the threshold for perfect matchings in random $r$-uniform hypergraphs.
However, since we have colored edges rather than $r$-tuples, it is not clear to us how to transfer our problem using Riordan's coupling.

The third problem is on the general transversal version of the problem.
If the host graphs are just random (hyper)graphs, then we are just taking i.i.d.~copies of random subgraphs of a complete graph (as the base graph) and consider transversal copies of our target subgraph.
However, a more general version allows the base (hyper)graph to be different while our Theorem~\ref{trans2} is only applicable when the host graphs are the same.
For instance, in a general transversal version of Theorem~\ref{treethm} (also Corollary~\ref{coro}), one can take $G_1, \dots, G_{n-1}$ as $n-1$ (not necessarily distinct and not necessarily the same) graphs each of which satisfies $\delta(G_i)\ge (1/2+\delta)n$, and consider a transversal spanning tree in the union of their sparsifications (that is, a spanning tree that contains exactly one edge from the sparsification of each $G_i$).
See~\cite{AnCh} for detailed discussions and results for the case of Hamilton cycles.

\section*{Acknowledgements}
The authors are indebted to Asaf Ferber and Huy Pham for their stimulating discussions at an earlier stage of this project.
The authors would like to thank the referees for comments that improve the presentation of the paper.

\appendix

\section{Well-definedness of the rainbow threshold}
\label{sec:mono}

\begin{lem} \label{mono}
Let $X$ be a finite set and let $X'=X \times [k]$, where $k$ is a positive integer. 
For $p\in (0,1)$, let $\mu_p'(S)$ be the discrete probability measure on $2^{X'}$ defined point-wise by
\[
\mu_p'(S): = \begin{cases}
(p/k)^t (1-p)^{|X|-t}  & \text{ if all $x_1, \dots, x_t$ are distinct,} \\
0 & \text{ otherwise,}
\end{cases} 
\]
for $S=\{(x_1, i_1), \dots, (x_t, i_t)\}\subseteq X'$.
Let $\F'\subseteq 2^{X'}$ be an increasing family, and assume $\F'\ne \emptyset, \ 2^{X'}$. Then $\mu_p'(\F')$ is strictly monotone increasing as $p$ increases.
\end{lem}

\begin{proof}
	We will show it for a more general measure which is no longer symmetric among the elements in $X$, and therefore we can prove it by considering the `local' contribution of each element in $X$. 
	
	Let $n = |X|$, and we may assume $X=[n]$. Let $p_1, p_2, \dots, p_n \in (0,1)$. 
	We define a probability measure $\mu_{\textbf{p}} = \mu_{p_1, p_2, \dots, p_n}$ point-wise by, for $S=\{(x_1, i_1), \dots, (x_t, i_t)\}\subseteq X'$,
\[
	\mu_{\textbf{p}} (S): = \begin{cases}
		\frac{1}{k^t} \prod_{i\in \{x_1, \dots, x_t\} } p_{i} \prod_{j\in X\setminus \{x_1, \dots, x_t\}} (1-p_j)
			& \text{ if all $x_1, \dots, x_t$ are distinct,} \\
		0 & \text{ otherwise,}
	\end{cases} 
\]
This generalizes our original model to include each vertex in randomly colored $X$ with its own probability, and $\mu_{\textbf{p}} = \mu_p'$ when $p_1 = p_2 = \cdots = p_n = p$.

Let $p_1<p_1'< 1$. Let $\mu$ denote $\mu_{\textbf{p}}$ and let $\mu^+$ denote $\mu_{p_1', p_2, \dots, p_n}$, i.e., replacing $p_1$ by $p_1'$ in the above formula. Next we will show that $\mu^+(\F')\ge \mu(\F')$.

Note that $\F'$ is increasing. For any $S \in \F'$ such that `$1$' never shows up as the first coordinate of the members of $S$, let $S_i := \{(1,i)\} \cup S$ for $i\in [k]$. Then $S_1, \dots, S_k$ are all in $\F'$ and they are distinct. That is, we can define a mapping $\phi : \F'\cap 2^{(X\setminus \{1\})\times [k]} \to \binom{\F'}k$, where $\phi (S) = \{S_1, \dots, S_k\}$. 

We claim that $\mu(\{S\}\cup \phi(S)) = \mu^+(\{S\}\cup \phi(S))$, for each $S\in \F'\cap 2^{(X\setminus \{1\})\times [k]}$. That is, as $p_1$ increases, the difference is canceled out within the tuple $(S, S_1,\dots, S_k)$. Indeed, for $S=\{(x_1, i_1), \dots, (x_t, i_t)\}\subseteq (X\setminus \{1\})\times [k]$,
\begin{align*}
		& \mu(\{S\}\cup \phi(S))  \\
	= &  \frac{1}{k^t} \prod_{i\in \{x_1, \dots, x_t\} } p_{i} \prod_{j\in X\setminus \{x_1, \dots, x_t\}} (1-p_j) 
		+ \sum_{l=1}^{k}\frac{1}{k^{t+1}} \prod_{i\in \{1, x_1, \dots, x_t\} } p_{i} \prod_{j\in X\setminus \{1, x_1, \dots, x_t\}} (1-p_j)  \\
	= & \left( \frac{1}{k^t} \prod_{i\in \{x_1, \dots, x_t\} } p_{i} \prod_{j\in X\setminus \{1, x_1, \dots, x_t\}} (1-p_j) \right) \left( (1-p_1) + p_1\right) \\
	= & \left( \frac{1}{k^t} \prod_{i\in \{x_1, \dots, x_t\} } p_{i} \prod_{j\in X\setminus \{1, x_1, \dots, x_t\}} (1-p_j) \right) \left( (1-p_1') + p_1'\right) \\
	= & \mu^+(\{S\}\cup \phi(S)).
\end{align*}

Clearly, $\phi (S)$ and $\phi (S')$ are disjoint subsets of $\F'$ for all distinct pairs $S,S'\in 2^{(X\setminus \{1\})\times [k]}$,
and are disjoint from $2^{(X\setminus \{1\})\times [k]}$. Let 
	\[
		\mathcal D := \F' \setminus \bigcup_{S\in \F'\cap 2^{(X\setminus \{1\})\times [k]}} (\{S\} \cup \phi(S)). \]
Then $\mathcal D$ (if not empty) and $\{S\} \cup \phi(S)$ for all $S \in \F'\cap 2^{(X\setminus \{1\})\times [k]}$ form a partition of $\F'$, 
and each $T\in \mathcal D$ has a member with `$1$' as its first coordinate. 
Hence, if $\mathcal D \ne \emptyset$, then for each $T = \{(1, i_1), (x_2, i_2), \dots, (x_t, i_t)\} \in \mathcal D$,
\begin{align} \label{Tdiff}
	\begin{split}
		 & \mu^+(T) - \mu(T) \\
		 = & \frac{1}{k^t} \cdot p_1'\prod_{i\in \{x_2, \dots, x_t\} } p_{i} \prod_{j\in X\setminus \{1, x_2, \dots, x_t\}} (1-p_j) 
		 	- \frac{1}{k^t} \prod_{i\in \{1, x_2, \dots, x_t\} } p_{i} \prod_{j\in X\setminus \{1, x_2, \dots, x_t\}} (1-p_j) \\
		= & \frac{1}{k^t} \cdot (p_1'-p_1) \prod_{i\in \{x_2, \dots, x_t\} } p_{i} \prod_{j\in X\setminus \{1, x_2, \dots, x_t\}} (1-p_j) > 0,
	\end{split}
\end{align}
as $p_i\in (0,1)$ for all $i$.

Now we are ready to look at the difference between $\mu^+$ and $\mu$ on the entire family $\F'$.
\begin{align} \label{mudiff}
	\begin{split}
		& \mu^+(\F') - \mu(\F') \\
	= & \sum_{S\in \F'\cap 2^{(X\setminus \{1\})\times [k]}} \left( \mu^+(\{S\}\cup \phi(S)) - \mu(\{S\}\cup \phi(S)) \right) 
		+  \sum_{T\in \mathcal D} \left( \mu^+(T) - \mu(T) \right) \\
	=  & \sum_{T\in \mathcal D} \left( \mu^+(T) - \mu(T) \right) \ge 0,
	\end{split}
\end{align}
and by \eqref{Tdiff}, the equation holds only if $\mathcal D =\emptyset$.

Thus by symmetry among $X$, \eqref{mudiff} shows that $\mu_{p_1, \dots, p_k} (\F)$ increases as $p_i$ increases for each $i\in [k]$; and hence, for any $0< p<p'< 1$, 
\[ \mu'_{p'}(\F') = \mu_{p', \dots, p'} (\F') \ge \mu_{p',  \dots, p', p} (\F') \ge \cdots \ge \mu_{p', p, \dots, p} (\F') \ge \mu_{p, \dots, p} (\F') = \mu'_p(\F'). \]
This shows the monotonity, and it is left to prove the strictness. 

Since $\F'\ne \emptyset$, it has a minimal element; and since $\F' \ne 2^{X'}$, $\emptyset$ is not its minimal element. Thus, without loss of generality, we may assume `$1$' is in a minimal element $M\in \F'$, 
then $M$ is in the set $\mathcal D$ at the iteration of increasing $p_1$ to $p_1'$. Thus by \eqref{Tdiff} and \eqref{mudiff}, 
\[ \mu'_{p'}(\F') \ge \mu_{p', p, \dots, p} (\F') > \mu_{p, \dots, p} (\F') = \mu'_p(\F').\]
Hence, 
$\mu_p'(\F')$ is strictly monotonic increasing as $p$ increases.
\end{proof}

\end{document}